\documentclass[leqno,12pt]{article} 
\usepackage{amsmath, amssymb}
\usepackage{amsthm} 
\theoremstyle{plain} 
\newtheorem{theorem}{\indent\sc Theorem}[section]

\newtheorem{corollary}[theorem]{\indent\sc Corollary}

\theoremstyle{definition} 

\makeatletter
\def\address#1#2{\begingroup
\noindent\parbox[t]{7.8cm}{%
\small{\scshape\ignorespaces#1}\par\vskip1ex
\noindent\small{\itshape E-mail address}%
\/: #2\par\vskip4ex}\hfill%
\endgroup}%
\makeatother
\title{\uppercase{Kernel theorem for the space of Beurling - Komatsu tempered
ultradistibutions}} 
\author{
\bigskip \\
\textsc{Z. Lozanov-Crvenkovi\'c and D. Peri\v{s}i\'{c}} 
}
\date{} 
\begin{document}

\maketitle

\footnote{ 
$^{*}$Thanks. }
\begin{abstract}
We give a simple proof of the Kernel theorem for the space of
tempered ultradistributions of Beurling - Komatsu type, using the
characterization of Fourier-Hermite coefficients  of the elements of
the space. We prove in details  that the test space of tempered
ultradistributions of Beurling - Komatsu type
 can be identified with the space of sequences of
ultrapolynomal falloff and its dual space  with the space of
sequences of ultrapolynomial growth. As a consequence of the Kernel
theorem we
 have that the Weyl transform can be  extended  on a space of
 tempered ultradistributions of Beurling - Komatsu type.
\end{abstract}

\section{Introduction}
Tempered ultradistributions
 spaces, as good spaces for harmonic analysis, have appeared in papers of many
 authors,
among others we mention
 Bj\"{o}rk \cite{BJ}, Wloka \cite{W},
Grudzinski \cite{GR},
  De Roever \cite{Roe},
 Kahspirovskij  \cite{KA} and Pilipovi\'c
 \cite{PT}, Matsuzawa, \cite{Matz},  Chung, \cite{Chung2}, \cite{Chung3}, Budin\v cevi\'c,
 Peri\v si\'c, Lozanov-Crvenkovi\'c \cite{DZM}  and Yoshino \cite{Yoshino}. These spaces have been
 developed
 in the framework of several ultradistribution theories: Komatsu's
 theory, Beurling-Bj\"{o}rk's theory, Cioranesku-Zsid\'o's and
 Braun-Meise-Taylor's theory.

In the paper we  study tempered ultradistributions following
Komatsu's approach to ultradistribution theory, introduced and
developed in his fundamental and inspiring papers \cite{K},
\cite{K2}, \cite{K3}.
 We study  the
space ${{\cal S}^{(M_p)}}'$ of tempered ultradistributions of
Beurling - Komatsu type,
   which  was introduced  and investigated in \cite{Varna} and \cite{Disertacija}. It
   is the generalization of the
space of tempered distributions ${\cal S}'$.

Beurling-Bj\"ork space ${\cal S}_{\omega}$, $\omega \in {\cal
M}_c$, introduced in \cite{BJ}, is equal to the space ${\cal
S}^{(M_p)}({\Bbb R}^d)$, where  $$M_p=
\sup_{\rho>0}\rho^pe^{-\omega(\rho)}.$$  The sequence satisfies
the conditions (M.1) and (M.3)', and it is in general different
from a Gevrey sequence. If we assume additionally that
$\omega(\rho) \geq C(\log \rho)^2$ for some $C>0$, then (M.2) is
satisfied.

The Pilipovi\'c
 space
$\Sigma_{\alpha}'$, (see \cite{PT}),  is an example of the space of
tempered ultradistributions, where $\{M_p\}$ is a Gevrey sequence $
p^{\alpha p},$ $p \in {\Bbb N}$, $\alpha >1$, which  appeared to be
interesting for the development of pseudodifferential operator
($\Psi$DO) theory in the framework of Time - Frequency analysis, see
the papers of Pilipovi\'c, Teofanov, \cite{PilTeo},
 \cite{PilTeo1}.

After introduction, in Section  2. we present without proof the
basic identification of the spaces ${\cal S}^{(M_p)}$ and ${{\cal
S}^{(M_p)}}'$ with the sequence spaces.
 In the sequence representation, the  space ${\cal
S}^{(M_p)}$ represents  the space of sequences of ultrapolynomial
falloff and ${{\cal S}^{(M_p)}}'$ represents the space of sequences
of  ultrapolynomial growth. We use these results to prove the Kernel
theorem for tempered ultradistributions.  As a consequence of the
Kernel theorem we
 have that the Weyl transform can be  extended  on a space of
 tempered ultradistributions of Beurling - Komatsu type and this
 gives the possibility to introduce the $\Psi$DO theory
 in the framework of tempered ultradistributions, where $\{M_p\}$
 is not necessarily a Gevrey sequence, but can be, for example
 $$
M_p = (p^{\nu}(\log  p )^{\mu})^p, \;\;\; p \in {\Bbb N},
$$
 where $\nu > 1$, and $\mu \in {\Bbb R}$, or $\nu = 1$ and $\mu
>1$.

In Section 3. we give the detailed proofs of characterizations of
Fourier-Hermite coefficients of elements of the spaces ${\cal
S}^{(M_p)}$ and ${{\cal S}^{(M_p)}}'$ which are only stated in
Section 2.
 This characterization appeared in several papers, \cite{Varna},
 \cite{Kaminski6}, \cite{PT}
 but in different forms and without a proof.
 Here we give  complete and detailed
proofs, based on the ideas from B. Simon's paper \cite{Simon} where
as a basic tool we use the harmonic oscillator wave functions -
Hermite functions.

\subsection{Notations}

Let $ \{M_{p},\;p \in  {\Bbb N}_{0}\}$ be a sequence of positive
numbers,  where $M_{0} =1$. An infinitely differentiable function
$\varphi$ on ${\Bbb R}^d$ is ultradifferentiable function of class
$(M_p)$ if on each compact set $K$ in ${\Bbb R}^d$ its derivatives
satisfy
$$
||\varphi^{(\alpha)}||_{C(K)} \leq C h^{|\alpha|}M_{|\alpha|},
\quad |\alpha|=0,1,...
$$
where we  use multi-index notation: for $d \in {\Bbb N}$, $\alpha
= (\alpha_{1},\alpha_{2},\dots,\alpha_{d}) \in {\Bbb N}_{0}^{d}$,
$x = (x_{1},x_{2},\dots,x_{d})$ $ \in {\Bbb R}^{d}$, we put
\begin{center}
 $|\alpha| =
\alpha_{1}+\alpha_{2}+\dots+\alpha_{d}$,
 $\;\;\;\alpha! = \alpha_{1}!\alpha_{2}!\cdots \alpha_{d}!$,
\\ $x^{\alpha} = x_{1}^{\alpha_{1}} x_{2}^{\alpha_{2}} \cdots
x_{d}^{\alpha_{d}}$, $|x|
=\sqrt{x_{1}^{2}+x_{2}^{2}+\cdots+x_{d}^{2}}$.
\end{center}
 For  $ x \in
{\Bbb R}^{d}$, $ \varphi(x) \in C^{\infty}({\Bbb R}^{d})$,
$$
\varphi^{(\alpha)}(x) = (\partial/\partial x)^{\alpha}\varphi(x) =
(\partial/\partial x_{1})^{\alpha_{1}}
  (\partial/\partial x_{2})^{\alpha_{2}}
  \cdots
  (\partial/\partial x_{d})^{\alpha_{d}}\varphi(x).
$$

 We impose the  conditions on the sequence $
\{M_{p},\;p \in  {\Bbb N}_{0}\}$ which are standard in
ultradistribution theory (for their detailed analysis see, for
example Komatsu's paper \cite{K}).
\begin{itemize}
\item[(M.1)]
$ M^{2}_{p} \leq M_{p-1} M_{p+1}, \quad  p = 1,2,\dots . $\\ ({\it
logarithmic convexity})
\\
\item[(M.2)]
There exist constants $A$, $H>0$ such that\\
 $ M_{p} \leq  A H^{p}
\min_{0 \leq q \leq p}M_{q}M_{p-q} ,
     \quad p = 0,1,\dots  \;
$\\
({\it stability under ultradifferential operators})
  \\
\item[(M.3)]
There exists  a constant $A>0$ such that \\
$
 \sum_{q=p+1}^{\infty } \dfrac{M_{q-1}}{M_{q}}
    \leq Ap\dfrac{M_{p}}{M_{p+1}},
    \quad p = 1,2,\dots .$
      \\({\it strong non-quasi-analyticity})
\end{itemize}

Some results remain valid, however, when (M.2) and (M.3) are
replaced by the following weaker conditions:

\begin{itemize}
\item[(M.2)']
There are constants $A$, $H>0$ such that\\
\\
$ M_{p+1} \leq A H^{p} M_{p},
   \quad p  = 0,1,\dots
 $  \\({\it stability under differential operators})
 \item[(M.3)']
$  B:= \sum^{\infty }_{p=1}
   \dfrac{M_{p-1}}{M_{p}} < \infty.$
 \\({\it non-quasi-analyticity})
\end{itemize}
We will often use the inequality
\begin{equation}\label{M*}
M_pM_q\leq M_0M_{p+q}, \quad p,q = 0,1,2,...,
\end{equation}
which follows from the condition (M.1), (see \cite[p. 45]{K}).

 The  associated function for the
        sequence   $\{M_{p},\; p \in {\Bbb N}_{0}\}$ is
 $$        M(\rho) =
        \sup_{p \in {\Bbb N}_0} \; \log
        \frac{\rho^{p}}{M_{p}}, \quad \rho > 0.
        $$
In the paper we will assume that the conditions (M.1), (M.2) and
(M.3)' are satisfied. The Gevrey sequence $ \{ p^{sp}\}_{p \in {\Bbb
N}_0}$, for $s >1,$ satisfies all the  above conditions and in that
case  $M(\rho) \sim \rho^{1/s}.$ The examples of such sequences are
also
\begin{equation}\label{sta}
M_p = (p^{\nu}(\log  p )^{\mu})^p, \;\;\; p \in {\Bbb N},
\end{equation}
where $\nu > 1$, and $\mu \in {\Bbb R}$, or $\nu = 1$ and $\mu
>1$. In this special case   $ M(\rho)=
\rho^{\frac{1}{2s}}(\log\,\rho)^{-\frac{t}{s}},  $  $\rho \gg 0$.
    There is a real need to consider that more general case. In
the study of the spaces od admissible data in the  Cauchy problems
it is necessary to consider ultradifferential classes wider than any
Gevrey  class. For example  Matsumoto  considered, see
\cite{Matsumoto} and references therein,  among others the case
(\ref{sta}).

Following Komatsu,
  we denote by  ${\cal
E}^{(M_p)}$ the space of all ultradifferentiable functions on ${\Bbb
R}^d$ of class $(M_p)$ and by ${\cal D}^{(M_p)}$ the subspace of
${\cal E}^{(M_p)}$, of all ultradifferentiable functions with
compact support.

If the sequence $ \{M_{p},\;p \in  {\Bbb N}_{0}\}$ satisfies
conditions (M.1) and (M.3)', by Denjoy-Carleman-Mandelbrojt
theorem ${\cal D}^{(M_p)}$ is not a trivial space. In this paper
we will always assume that   these conditions are satisfied.

We denote by ${{\cal D}^{(M_p)}}'$ the strong dual space of ${\cal
D}^{(M_p)}$, and call its elements  ultradistributions of Beurling -
Komatsu type.

\subsection{Space of tempered ultradistributions of Beurling -
Komatsu type}

The space ${{\cal S}^{(M_p)}}'$  of tempered ultradistributions of
type $(M_p)$,  is a subspace of the space of  Beurling-Komatsu
ultradistributions. The ${\cal D}^{(M_p)}$ is dense in ${\cal
S}^{(M_p)}$ and  set ${\cal S}^{(M_p)} \backslash {\cal D}^{(M_p)}$
is nonempty. Moreover: In the following diagram  the arrows denote
continuous inclusions:
$$
\begin{array}{ccccc}
{\cal D}' & \longleftarrow& {\cal S}' & \longleftarrow& {\cal
E}'\\
      \downarrow      &   &      \downarrow         &  &\downarrow \\
 {{\cal D}^{(M_p)}}' & \longleftarrow & {{\cal S}^{(M_p)}}'
&\longleftarrow&{{\cal
E}^{(M_p)}}'\\
    &   &      \downarrow         &  & \\
      &   &     {\cal G}'      &  & \\
\end{array}
$$
where ${\cal S}'$ is the Schwartz space tempered distributions and
${\cal G}'$ is the space  of extended Fourier hyperfunctions
(defined as in \cite{Chung2}). For the proof see \cite{Prim}.

Let us now give the precise definitions.

By ${\cal S}^{M_p,m}, $  $m >0$, we denote the Banach space of
smooth functions $\varphi $ on ${\Bbb R}^{d}$, such that for some
constant $C>0$
\begin{equation}\label{3.1}
||x^{\beta}\varphi^{(\alpha)}||_2 \leq C
\frac{M_{|\alpha|}M_{|\beta|}}{m^{|\alpha|+|\beta|}},\;\;
\mbox{for every $\alpha, \beta \in {\Bbb N}_{0}^{d}$},
\end{equation}
where $||\; ||_2$ is the standard norm in space $L^2({\Bbb
R}^d))$, equipped with the norm
\begin{equation}\label{M11}
s_{m}(\varphi) = \sup_{\alpha,\beta \in {\Bbb N}_{0}^{d}}
\frac{m^{|\alpha|+|\beta|}}{M_{|\alpha|}M_{|\beta|}}
||x^{\beta}\varphi^{(\alpha)}(x)||_2.
\end{equation}
The test space for the space of tempered ultradistributions is
defined by
$$
{\cal S}^{(M_{p})} = \operatornamewithlimits{proj\;lim}_{m
\rightarrow \infty} {\cal S}^{M_p,m}.
$$
It is a Frechet space.

The strong dual ${{\cal S}^{(M_{p})}}'$ of the space ${\cal
S}^{(M_{p})}$ is the space of tempered ultradistributions  of
Beurlinger - Komatsu type. The space of tempered ultradistributions
is a good space for Harmonic analysis, since the Fourier transform
is an isomorphism of ${{\cal S}^{(M_{p})}}'$ onto itself, (see
\cite{Kaminski7}). It is defined on ${\cal S}^{(M_{p})}$  by
$$
{\cal F}\varphi(\xi)= \int_{{\Bbb R}} e^{-2\pi i x \xi }f(x)dx,
               \quad \varphi \in {\cal S}^{(M_{p})},
$$
and on ${{\cal S}^{(M_{p})}}'$ by
$$
 \langle {\cal F}f,\varphi \rangle  =  \langle f,{\cal F}\varphi \rangle ,
$$
for $f \in {{\cal S}^{(M_{p})}}'$ and $ \varphi\in  {\cal
S}^{(M_{p})}$.

It has a lot of nice properties, for example, (see \cite{Kaminski7}
and \cite{Disertacija} for more details), if
 $$
 P(x,D)=\sum a_{\mu,\nu}(-1)^{\nu}D^{\nu}x^{\mu}
 $$
 where $a_{\mu,\nu}$ are complex numbers such that there exists $L>0$ and $C>0$ and
 $$
 |a_{\mu,\nu}|\leq S \frac{L^{\mu+\nu}}{M_{\mu}M_{\nu}}, \quad \mu.\nu \in {\Bbb N},
 $$
 then
 $$
 {\cal F}(P(\cdot,D)f(\cdot))(\xi) = P(-D,\xi)({\cal F}f)(\xi), \quad \xi \in {\Bbb R},
 $$
 for each $f \in {{\cal S}^{(M_{p})}}'$.

\section{Main Results}

Komatsu proved, in his papers, \cite{K2} and \cite{K3} the Kernel
theorem for
 ultradistributions and ultradistributions with compact support
  which are an analogue of L. Schwartz Kernel theorem for
distributions.

The main result of the paper is the Kernel theorem for tempered
ultradistributions. Applying the theorem we show that the Weyl
transform, which is a representation of the convolution algebra on
$L^1$ functions over the reduced Heisenberg group $H_d^{red}$,
extends uniquely to a bijection from ${{\cal S}^{(M_{p})}}'({\Bbb
R}^{2l})$ to the space of continuous linear maps from ${{\cal
S}^{(M_{p})}({\Bbb R}^{l})}$ to
 ${{\cal S}^{(M_{p})}}'({\Bbb R}^{l})$.

The simple proof of the Kernel theorem for tempered
ultradistributions, which will be given in this section, depend on
the characterization of the Fourier-Hermite coefficients of elements
of ${{\cal S}^{(M_{p})}}$ and ${{\cal S}^{(M_{p})}}'$. In this
section we will state the theorems which give these
characterizations and prove the Kernel theorem. We follow the ideas
from the paper \cite{Simon} of B. Simon.

The elements of the space ${\cal S}^{(M_p)}$ are  in $L^2$, thus
they have in the space  the Hermite  expansion $\sum a_n(\varphi)
h_n$ where $h_n$ are Hermite functions:
$$
h_n(x) = (-1)^n \pi^{-\frac{1}{4}}
2^{-\frac{1}{2}n}(n!)^{-\frac{1}{2}}e^{\frac{1}{2}x^2}\Big(
\frac{d}{dx} \Big)e^{-x^2}, \quad n \in {\Bbb N},
$$
which are    harmonic oscillator wave functions, since
$$
\Big(-\displaystyle\frac{d}{dx^2} + x^2 \Big)h_n = (2n+1)h_n.
$$

 Let $\varphi \in {\cal
S}^{(M_p)}$, the Hermite coefficients of $\varphi$ are
$$a_n(\varphi) =\int_{\Bbb R} \varphi(x) h_n(x)dx, \quad n \in {\Bbb N}_0.$$
The sequence of the Fourier-Hermite coefficients $\{a_n\}_{n\in
{\Bbb N}_0}$ of $\varphi$ we call the {\bf Hermite representation
of} $\varphi$.

The first main result of the paper is the Hermite representation of
the space
 $ {\cal S}^{(M_p)}({\Bbb R})$:
\begin{theorem}\label{grejalica}
(a) If $\varphi \in  {\cal S}^{(M_p)}({\Bbb R})$, then for every
$\theta
>0$
 \begin{equation}
 ||\varphi||_{\theta}= \sum_{n =0}^{\infty} |a_n(\varphi)|^2
 \exp[2M(\theta \sqrt{2n+1})] < \infty.
 \end{equation}
(b) Conversely, if a sequence $\{a_n\}_{n \in {\Bbb N}_0 }$ satisfy
that for every $\theta >0$,
 $$ \sum_{n =0}^{\infty} |a_n|^2 \exp[2M(\theta \sqrt{2n+1})] <
 \infty,$$
   then the series $ \sum a_n h_n  $  converges in
 the space $ {\cal S}^{(M_p)}({\Bbb R})$ to some $\varphi$,
 and the sequence $\{a_n\}_{n \in {\Bbb N}_0 }$
 is the Hermite representation of $\varphi$.
\end{theorem}

From the proof of  Theorem \ref{grejalica}, which will be given in
the next section, it follows that that the test space of tempered
ultradistributions of Beurling - Komatsu type ${\cal
S}^{(M_p)}({\Bbb R})$ can be identified with the space $s^{(M_p)}$
of sequences of ultrapolynomal falloff, where we define it by
$$
s^{(M_p)}=\{\{a_n\}_{n \in {\Bbb N}_0}\,|\, \forall \theta >0, \,
 \sum_{n =0}^{\infty} |a_n|^2 \exp[2M(\theta
\sqrt{2n+1})] <
 \infty \},
$$
and equip it  with the topology induced by the family  of
seminorms
$$\{
||\{a_n\}||_{\theta} =\sum_{n =0}^{\infty} |a_n|^2 \exp[2M(\theta
\sqrt{2n+1})],\;\; \theta >0 \}.$$
 More precisely, the map
$$
{\cal S}^{(M_p)}({\Bbb R}) \rightarrow s^{(M_p)}, \quad \varphi
\mapsto \{a_n(\varphi)\}_{n \in {\Bbb N}_0}, \quad a_n(\varphi) =
 \langle \varphi, h_n \rangle
$$
is a topological isomporphism. It is easy to prove that the space
$s^{(M_p)}$ is  nuclear and therefore space ${\cal
S}^{(M_p)}({\Bbb R})$ is nuclear.

If $f \in {{\cal S}^{(M_p)}}'({\Bbb R}) $, then  $a_n(f) =
 \langle f,h_n \rangle $, $n = 0,1,2,...$, are the Fourier-Hermite
coefficients of $f$, and the sequence $\{a_n(f)\}_{n \in {\Bbb
N}_0}$, is  the {\bf Hermit representation} of tempered
ultradistribution $f$.

Hermite representation of the space
 $ {{\cal S}^{(M_p)}}'({\Bbb R})$ is given by:
\begin{theorem}\label{ljilja}

 (a) Let
 $f \in {{\cal S}^{(M_p)}}'({\Bbb R})$. Then
 for some
$\theta >0$
\begin{equation}\label{zaga}
\sum_{n=0}^{\infty}|a_n(f)|^2\exp[-2M(\theta\sqrt{2n+1})] < \infty
\end{equation}
and $ \langle f,\varphi \rangle  = \sum_{n=0}^{\infty} a_n(f)
a_n(\varphi)$, where the sequence $\{a_n(\varphi)\}_{n \in {\Bbb
N}_0}$ is the Hermite representative of $\varphi \in {\cal
S}^{(M_p)}({\Bbb R})$.

(b) Conversely, if a sequence $\{b_n\}_{n \in {\Bbb N}_0}$ satisfy
$$
\sum_{n =0}^{\infty}|b_n|^2\exp[-2M(\theta\sqrt{2n+1})]<\infty
\quad \text{for some}\quad \theta >0,
$$
then the sequence $\{b_n\}_{n \in {\Bbb N}_0}$ is the Hermite
representation of some $f \in  {{\cal S}^{(M_p)}}'({\Bbb R})$ and
 Parseval equation holds
 $ \langle f,\varphi \rangle  =\sum_{n=0}^{\infty}a_n(f)a_n(\varphi).$
 \end{theorem}

Although we have stated the results for ${\cal S}^{(M_p)}({\Bbb
R})$ and ${{\cal S}^{(M_p)}}'({\Bbb R})$ in  one dimensional case,
analogous results hold in the multidimensional case, for ${\cal
S}^{(M_p)}({\Bbb R}^l)$ and ${{\cal S}^{(M_p)}}'({\Bbb R}^l)$. The
Hermite functions in multidimensional case are defined by
$$
h_n(x)= h_{n_1}(x_1)h_{n_2}(x_2)\cdots h_{n_d}(x_d), \quad x
=(x_1,x_2,...x_d)\in {\Bbb R}^d,
$$
where $n=(n_1,n_2,...,n_d) \in {\Bbb N}^d$. By $h_{(n,k)}$ we denote
$$
h_{(n,k)}=h_{n_1}(x_1)h_{n_2}(x_2)\cdots h_{n_l}(x_l)
h_{k_1}(x_{l+1})h_{k_2}(x_{l+2})\cdots h_{k_s}(x_{l+s}),\,\,
(n,k)\in {\Bbb N}^l_0\times {\Bbb N}^s_0.
$$

\subsection{The Kernel theorem}

The Kernel theorem for tempered ultradistributions states
essentially that
 every continuous linear map ${\cal K}$ of the space $({\cal S}^{(M_p)}({\Bbb R}^l))_x$
 of test functions in some variable $x$, into the space $({{\cal S}^{(M_p)}}'({\Bbb R}^s))_y$
of tempered ultradistributions in a second variable $y$, is given
by a unique tempered ultradistributions  $K$ in both variables.

\begin{theorem}[{\bf Kernel theorem}] Every jointly continuous bilinear
 functional $K$ on
  $ {{\cal S}^{(M_p)}}({\Bbb R}^l)\times {{\cal
 S}^{(M_p)}}({\Bbb R}^s)$ defines a linear map ${\cal K} : {\cal
 S}^{(M_p)}({\Bbb R}^s)\rightarrow {{\cal
 S}^{(M_p)}}'({\Bbb R}^l)$ by
\begin{equation}\label{kraj}
 \langle  {\cal K}\varphi, \psi \rangle  = K(\psi \otimes\varphi),\quad
\text{for}\;\;\,\varphi \in {\cal S}^{(M_p)}({\Bbb R}^s), \psi\in
{\cal S}^{(M_p)}({\Bbb R}^l).
\end{equation}
and $(\varphi\otimes\psi)(x,y)=\varphi(x)\psi(y)$, which is
continuous in the sense that ${\cal K} \varphi_j \rightarrow 0 $ in
${{\cal
 S}^{(M_p)}}'({\Bbb R}^l)$ if $ \varphi_j \rightarrow 0 $ in ${{\cal
 S}^{(M_p)}}({\Bbb R}^s)$.

 Conversely, for linear map ${\cal K} : {\cal
 S}^{(M_p)}({\Bbb R}^s)\rightarrow {{\cal
 S}^{(M_p)}}'({\Bbb R}^l)$ there is unique
 tempered ultradistribution $K \in {{\cal S}^{(M_p)}}'({\Bbb R}^{l+s})$ such that (\ref{kraj}) is valid.
 The tempered ultradistribution $K $ is called the kernel of ${\cal K}$.
 \end{theorem}
Proof.  If $K$ is a jointly continuous bilinear functional on $
{{\cal S}^{(M_p)}}({\Bbb R}^l)\times {{\cal
 S}^{(M_p)}}({\Bbb R}^s)$, then (\ref{kraj}) defines a tempered
 ultradistribution $({\cal K}\varphi) \in {{\cal
 S}^{(M_p)}}'({\Bbb R}^l)$ since
 $
 \psi \mapsto K(\psi\otimes\varphi)
 $
is continuous. The mapping ${\cal K}:{{\cal S}^{(M_p)}}({\Bbb
R}^{s})\rightarrow {{\cal S}^{(M_p)}}'({\Bbb R}^{l})$ is
continuous since the mapping $ \varphi \mapsto
K(\psi\otimes\varphi)
 $
is continuous.

Let us prove the converse. To prove the existence we define a
bilinear form $B$ on ${{\cal S}^{(M_p)}}'({\Bbb R}^{l})\otimes
{{\cal S}^{(M_p)}}'({\Bbb R}^{s})$ by
$$
B(\varphi,\psi) =  \langle {\cal K}\psi,\phi \rangle , \quad
\psi\in {\cal S}^{(M_p)}({\Bbb R}^l),\; \varphi \in {\cal
S}^{(M_p)}({\Bbb R}^s).
$$
The form $B$ is a separately continuous  bilinear form on the
product ${\cal S}^{(M_p)}({\Bbb R}^l)\times {\cal S}^{(M_p)}({\Bbb
R}^s)$ of Frechet spaces and therefore it is jointly continuous,
see \cite{Trev}.

 Let $C>0$, $\theta \in {\Bbb R}^l_+$, $\nu \in {\Bbb R}^s_+$ be chosen so that
\begin{equation}\label{B1}
|B(\varphi,\psi)|\leq C||\varphi||_{\theta}||\psi||_{\nu},
\end{equation}
and let
\begin{equation*}
t_{(n,k)}= B(h_n,h_k), \quad n \in {\Bbb N}^l,k\in {\Bbb N}^s.
\end{equation*}
Since $B$ is jointly continuous on ${{\cal S}^{(M_p)}}({\Bbb
R}^l)\times {{\cal
 S}^{(M_p)}}({\Bbb R}^s)$, for
 $\varphi=\sum a_n h_n$ and
$\psi=\sum b_k h_k$ we have that
\begin{equation*}
B(\varphi,\psi)=\sum t_{(n,k)}a_n b_k.
\end{equation*}
On the other hand, for $ (n,k) \in {\Bbb N}^l \times {\Bbb N}^s$ and
$(\theta,\nu)\in {\Bbb R}^l \times{\Bbb R}^s$, by (\ref{B1}) we have
\begin{equation*}
\begin{split}
 |t_{(n,k)}|& \leq
C||h_{n}||_{\theta}||h_k||_{\nu}=\\
&= ||h_{n_1}||_{\theta_1} ||h_{n_2}||_{\theta_2}\cdots
||h_{n_l}||_{\theta_l}\; ||h_{k_1}||_{\nu_1} ||h_{k_2}||_{\nu_2}
\cdots||h_{k_s}||_{\nu_s}=\\
&=\exp[2\sum_{i=1}^{l}M(\theta_i\sqrt{n_i})]\exp[2\sum_{j=1}^{s}M(\nu_j\sqrt{k_j})].
\end{split}
\end{equation*}

Thus, from Theorem \ref{ljilja}  it follows that  the sequence
$\{t_{(n,k)}\}_{(n,k)}$ is a Hermite representation of a tempered
ultradistribution $K\in{{\cal S}^{(M_p)}}'({\Bbb R}^l\times {\Bbb
R}^s)$. Thus
\begin{equation}\label{dosta}
 \langle K,\varphi \rangle = \sum t_{(n,k)}c_{(n,k)},
\end{equation}
for  $\varphi=\sum c_{(n,k)}h_{n,k}\in {\cal S}^{(M_p)}({\Bbb
R}^{l+s})$.

If $\varphi=\sum a_n h_n \in {\cal S}^{(M_p)}({\Bbb R}^{l})$ and
$\psi=\sum b_k h_k \in {\cal S}^{(M_p)}({\Bbb R}^{s})$ then
$\varphi\otimes\psi$ has the Hermite representation
$\{a_nb_k\}_{(n,k)}$ and we have that for tempered
ultradistribution $K$ defined by (\ref{dosta})
$$
K(\varphi\otimes\psi)=\sum_{(n,k)}t_{(n,k)}a_nb_k = B(\varphi,\psi),
$$
so $K=B$.  This proves the existence.

The uniqueness follows from the fact that $K$ is completely
determined by its Hermite representation $\{ \langle K,h_{(n,k)}
\rangle \}_{(n,k)}$ and the fact that for every $(n,k)\in {\Bbb
N}^l \times {\Bbb N}^s $
$$
 \langle K,h_{(n,k)} \rangle  =  \langle K,h_n\otimes h_k \rangle = B(h_n,h_k) =
t_{(n,k)}.
$$
QED

The unitary representation $\rho$ of the reduced Heisenberg group
$H_n^{red}$ might be considered as a representation of $L^1({\Bbb
R}^{2l})$, with a nonstandard convolution structure. Accordingly,
for $F \in L^1({\Bbb R}^{2l})$ let us define (see \cite{Foland}) a
bounded operator $\rho(F)$  on $L^2({\Bbb R}^l)$ (it is sometimes
called the Weyl transform of $F$) as
\begin{equation}\label{ro}
\rho(F) \varphi(x)= \int\int F(y-x,q)e^{\pi i q(x+y)}
\varphi(y)dydq.
\end{equation}
In other words $\rho(F)$ is an integral operator with kernel
$$
K_F(x,y)= {\cal F}_2^{-1}F(y-x,\frac{y+x}{2}),
$$
where ${\cal F}_2$ denotes Fourier transform in the second variable
and ${\cal F}_2^{-1}$ its inverse transform.  Calling $\rho(F)$ the
Weyl transform is historically inaccurate. In fact,
 the Weyl transform of $F$ should be $\rho({\cal F}(F))$.
 \begin{corollary}
 The map $\rho$ from $L^1({\Bbb R}^{2l})$ to the space of bounded
 operators on $L^2({\Bbb R}^l)$, defined by (\ref{ro}), extends
 uniquely to a bijection from ${{\cal S}^{(M_{p})}}'({\Bbb
 R}^{2l})$to the space of continuous linear
maps from ${{\cal S}^{(M_{p})}({\Bbb R}^{l})}$ to
 ${{\cal S}^{(M_{p})}}'({\Bbb R}^{l})$.
 \end{corollary}

\begin{proof}
The kernel $K_F$ is well defined when $F$ is an arbitrary
ultradisrtibution  and belongs to the space ${{\cal S}^{(M_p)}}'$.
So the Weyl transform extends to a bijection from ${{\cal
S}^{(M_{p})}}'({\Bbb R}^{2l})$ to the space of continuous linear
maps from ${{\cal S}^{(M_{p})}({\Bbb R}^{l})}$ to
 ${{\cal S}^{(M_{p})}}'({\Bbb R}^{l})$.

  The uniqueness of the extension follows from the kernel theorem
 since every continuous linear map form ${{\cal S}^{(M_{p})}({\Bbb
 R}^{l})}$  to ${{\cal S}^{(M_{p})}}'({\Bbb R}^{l})$ is of the
 form (\ref{kraj}).
\end{proof}

\section{Proofs of Theorems \ref{grejalica} and \ref{ljilja}}

 In the proof of Theorem \ref{grejalica} we will use the following facts:

 (i) Since the sequence $m_n = M_n/M_{n-1}$, $n = 1,2,...$ is
 increasing (which is equivalent to the condition (M.1), see \cite[p.
 50]{K}), we have
\begin{equation}\label{GG}
 \frac{k}{m_k}\leq \sum_{k=0}^{\infty}\frac{1}{m_k}=:B < \infty,
\end{equation}
 and
\begin{equation}\label{SS}
 \frac{k!}{n!}\frac{M_n}{M_k}=\frac{k}{m_k}\frac{k-1}{m_{k-1}}\cdots\frac{n+1}{m_{n+1}}\leq
 B^{k-n}, \quad n,k \in {\Bbb N}, n\leq k.
\end{equation}

From above and Stirling formula
\begin{equation}\label{Stirling}
 n^n = n! \frac{1}{\sqrt{2\pi n}}e^ne^{-g(n)}\leq Ce^n n!,\quad |g(n)|\leq
 \frac{1}{12n},
\end{equation}
 it follows that
\begin{equation}\label{DD}
 \frac{k^k}{n^n}\frac{M_n}{M_k}\leq C e^k\sqrt{n}B^{k-n}, \quad
 k,n \in {\Bbb N}, n\leq k.
\end{equation}

(ii)   For  $n,m \in {\Bbb N}$
\begin{equation}\label{2.5}
x^mh_n(x) = 2^{-m/2}\sum_{k=0}^m \alpha^{(n)}_{k,m}h_{n-m+2k}(x),
\quad x \in \Bbb R,
\end{equation}
where
\begin{equation}\label{2.6}
0 \leq |\alpha^{(n)}_{k,m}|\leq {m \choose k}
((2n+1)^{m/2}+m^{m/2})
\end{equation}

(iii) If we denote ${\cal R}=\Big(-\displaystyle\frac{d}{dx^2} +
x^2 \Big) $, then for $\varphi \in {\cal S}$ and $k \in \Bbb N $
it holds
\begin{equation}\label{ocena}
{\cal R}^k \varphi=\sum_{n=0}^{k}\sum_{p=0}^{2n}
C^{(k)}_{p,2n-p}\,x^p\varphi^{(2n-p)}, \quad |C^{(k)}_{p,2n-p}|\leq
10^k k^{k-n}.
\end{equation}
The (ii) and (iii) were proved by mathematical induction in
\cite{KA}.

\begin{proof}  (Theorem 2.1) We give a detailed proof  in one dimensional
case.

a) Let us first prove that  if $\varphi \in  {\cal S}^{(M_p)}$,
then for every $\theta >0$
 \begin{equation}\label{teta norma}
  \sum_{n=0}^{\infty} |a_n|^2
 \exp[2M(\theta \sqrt{2n+1})] < \infty.
 \end{equation}
 Let
$\varphi \in {\cal S}^{(M_p)}$. From (\ref{ocena}), (\ref{3.1})
(M.1), (\ref{M*}), Stirling formula, (M.2),  (\ref{SS}), (\ref{DD})
it follows that there exists $C$ such that for each $m
> 0$ and each $k \in {\Bbb N}_0$,

\begin{equation*}
\begin{split}
\Big(\sum_{n \in {\Bbb N}_0}&|a_n|^2(2n+1)^{2k}\Big)^{1/2}=||{\cal
R}^k
\varphi|| _2\leq \sum_{n=0}^k \sum_{p=0}^{2n}C^{(k)}_{p,2n-p}||x^p\varphi^{(2n-p)}||_2\\
&\leq C  \sum_{n=0}^k \sum_{p=0}^{2n} 10^k k^{k-n}m^{-2n}
M_pM_{2n-p}\\
& \leq C  \sum_{n=0}^k \sum_{p=0}^{2n}10^k\frac{k^k}{n^n}
        m^{-2n}M_{2n}\\
&\leq  C  \sum_{n=0}^k \sum_{p=0}^{2n}
10^k\frac{k!}{n!}e^{k-n+\frac{1}{12n}}\sqrt{n}\,m^{-2n}M_{2n}\\
&\leq C\,20^k e^k (1+H^2)^k m^{2k}M_{2k}  \sum_{n=0}^k
(2n+1)\frac{k!}{n!}\frac{M_n}{M_k}\frac{M_n}{M_k}\frac{1}{m^{2(k-n)}}\\
&\leq C\,80^k e^k (1+H^2)^k m^{2k}M_{2k}
\sum_{n=0}^k\frac{1}{2^n}B^{k-n}\frac{M_n}{M_k}\left(\frac{1}{m^2}\right)^{k-n}\\
& \leq C\,80^k e^k (1+H^2)^k\,B^k(1+\frac{1}{B})^k m^{2k}
M_{2k}\times\\
& \qquad\qquad\qquad\times\left(
\sum_{n=0}^{[\frac{1}{m^2}]}\frac{H^{k+1}}{2^n}\frac{M_n}{M_{k+1}}
\frac{(k+1)^{k+1}}{n^n}
+\sum_{n=[\frac{1}{m^2}]+1}^{k}\frac{H^{n+k+1}}{2^n}\frac{M_{n-1}}{M_{k+1}}\frac{(k+1)^{k+1}}{(n-1)^{n-1}}\right)\\
&\leq C\,80^k e^k (1+H^2)^k B^k (1+\frac{1}{B})^k m^{2k}
M_{2k}\times\\
&  \qquad\qquad\qquad\times\left(
\sum_{n=0}^{[\frac{1}{m^2}]}\frac{H^{k+1}}{2^k}e^{k+1}\sqrt{n}\,B^{k+1-n}
+ \sum_{n=[\frac{1}{m^2}]+1}^{k}\frac{H^{n+k+1}}{2^n}e^{k+1}\sqrt{n-1}\,B^{k+1-n}\right)\\
& \leq C \,160^ke^{2k}(1+H^2)^k (1+H)^{2k}B^{2k}(1+\frac{1}{B})^{2k}m^{2k}M_{2k} \\
&\leq C\, m_{1}^{2k} M_{2k}
\end{split}
\end{equation*}
where $m_1 = \sqrt{160(1+H^2)}\,e(1+H)B(1+\frac{1}{B}) m$.

 It follows that, for $2k = \alpha +2$
$$
|a_n|^2(2n+1)^{\alpha+2}\leq C m_1^{\alpha+2}M^2_{\alpha+2}.
$$
Applying (M.2) we obtain
$$
|a_n|^2(2n+1)^{\alpha+2}\leq C m_1^{\alpha+2}A^2 H^{2(\alpha +2)}
M^2_{\alpha}
$$
which implies that
 for each $\alpha$, $n \in {\Bbb
N}_0$ and $\theta = m_1 H$, there exists $C$, such that
$$
\frac{|a_n|^2 \theta^{2\alpha}(2n+1)^{\alpha}}{M_{\alpha}^2}\leq
\frac{C}{(2n+1)^2},
$$
which implies
$$
|a_n|^2 \exp[2M(\theta\sqrt{2n+1})]=|a_n|^2 \sup_{\alpha \in {\Bbb
N}_0}\frac{\theta^{2\alpha}(2n+1)^{\alpha}}{M^2_{\alpha}} \leq
\frac{C}{(2n+1)^2}.
$$
Therefore
$$ \sum_{n \in {\Bbb N}_0} |a_n|^2
\exp[2M(\theta \sqrt{2n+1})] < \infty,
$$
for every $\theta >0$.

b) Let us now prove that if
\begin{equation}\label{mob}
 \sum_{n
\in {\Bbb
 N}_0} |a_n|^2 \exp[2M(\theta \sqrt{2n+1})] < \infty,
\end{equation}
 for every $\theta >0$,   then  the series $ \sum a_n h_n  $ converges in
 the space $ {\cal S}^{(M_p)}$.

Suppose that  inequality (\ref{mob}) holds for every  $\theta >0$.
Applying (\ref{2.5}), (\ref{2.6}),  Cauchy-Schwartz inequality  and
\begin{equation}\label{kvadrat}
\exp[-M(\theta\sqrt{2n+1})]\leq \frac{M_4}{\theta^4 (2n+1)^2},
\end{equation}
 which follows from the definition of the associated function.
 We get that for every $\theta>0$ there exist a
constant $C$ such that for each $p \in {\Bbb N}_0$
\begin{equation*}
\begin{split}
&||x^p\varphi||_2 \leq 2^{-p/2}||\sum_{n \in {\Bbb
N}_0}a_n\Big(\sum_{k \leq
p}\alpha^{(n)}_{k,p}h_{n-p+2k}\Big)||_2 \\
&\leq 2^{-p/2}\sum_{n \in {\Bbb N}_0}a_n
((2n+1)^{p/2}+p^{p/2})\Big(\sum_{k \leq p} {p \choose k}
\Big)||h_{n-p+2k}||_2\\
&\leq 2^{p/2}\sum_{n \in {\Bbb N}_0}a_n
\exp[M(\theta\sqrt{2n+1})-M(\theta\sqrt{2n+1})]
 ((2n+1)^{p/2}+p^{p/2})\\
&\leq 2^{p/2}\Big(\sum_{n \in {\Bbb N}_0}|a_n|^2
\exp[2M(\theta\sqrt{2n+1})]\Big)^{1/2}\times\\
& \qquad \qquad \qquad \times \Big(\sum_{n \in {\Bbb
N}_0}\exp[-2M(\theta\sqrt{2n+1})]((2n+1)^{p/2}+p^{p/2})^2
\Big)^{1/2}\\
&\leq C 2^{p/2}\Big(\sum_{n \in {\Bbb
N}_0}\exp[-2M(\theta\sqrt{2n+1})]((2n+1)^{p/2}+p^{p/2})^2
\Big)^{1/2}\\
&\leq C 2^{p/2}\sup_{n \in {\Bbb
N}_0}\Big(((2n+1)^{p/2}+p^{p/2})\exp[-\frac{1}{2}M(\theta\sqrt{2n+1})]\Big)\times\\
& \qquad \qquad \qquad \times\Big(\sum_{n \in {\Bbb
N}_0}\exp[-M(\theta\sqrt{2n+1})] \Big)^{1/2}\\
&\leq C\left(\Big(\frac{2}{\theta}\Big)^{p/2}M_p\Big(\theta^{p/2}
\sup_{n}\frac{(2n+1)^{p/2}\exp[-\frac{1}{2}M(\theta\sqrt{2n+1})]}{M_p}\Big)+\frac{\theta^{p/2}p^{p/2}}{M_p}\right)\times\\
& \qquad \qquad \qquad \times \frac{M_4}{\theta^4}\sum_{n\in {\Bbb
N}_0}\frac{1}{(2n+1)^2}
\end{split}
\end{equation*}
Since
$$
\frac{\theta^{p/2}p^{p/2}}{M_p} \leq
\frac{p!e^p\theta^{p/2}}{M_p}\rightarrow 0,\quad \text{as} \quad
p\rightarrow\infty,
$$
and
\begin{equation*}
\begin{split}
\sup_{n \in {\Bbb
N}_0}\frac{\theta^{p/2}(2n+1)^{p/2}\exp[-\frac{1}{2}M(\theta\sqrt{2n+1})]}{M_p}
&=\frac{1}{M_p}\sup_{n \in {\Bbb N}_0}\Big(
\frac{\theta^{p/2}(2n+1)^p}{\exp[M(\theta\sqrt{2n+1})]}\Big)^{1/2}\\
&\leq \frac{\sqrt{M_p}}{M_p}\rightarrow 0,\quad \text{as} \quad
p\rightarrow\infty,
\end{split}
\end{equation*}
which follows from \cite[(3.3)]{K}, we have that for each  $\theta
>0$  there exists $C>0$ such that
\begin{equation}\label{2.14}
||x^p\varphi||_2\leq C\Big(\sqrt{2/\theta}\Big)^p M_p,\quad p \in
{\Bbb N}_0.
\end{equation}
From above and properties of
 the Fourier transform we obtain that for each  $\theta
>0$  there exists $C>0$ such that
\begin{equation}\label{2.15}
||\varphi^{(q)}||_2
 =\frac{1}{\sqrt{2\pi}}||x^q{\cal F}(\varphi)||_2=\frac{1}{\sqrt{2\pi}}||x^q
 \sum_{k \in {\Bbb N}_0}a_k
h_k||_2=\frac{1}{\sqrt{2\pi}}||x^q\varphi||_2
 \leq
C\Big(\sqrt{2/\theta}\Big)^q M_q.
\end{equation}

If $p,q \in {\Bbb N}_0$ and $\gamma = \min(q, 2p)$, by using
(\ref{2.14}),(\ref{2.15}), (M.1), (M.3)' and (M.2) we get
\begin{equation}\label{2.17}
\begin{split}
&\Big(||x^{p}\varphi^{(q)}||_2\Big)^2=(x^{p}\varphi^{(q)},x^{p}\varphi^{(q)})_{L^2}
=|((x^{2p}\varphi^{(q)})^{(q)},\varphi)_{L^2}|\\
&\leq \Big|\sum_{k=0}^{q}{q\choose
k}\frac{(2p)!}{(2p-k)!}(x^{2p-k}
\varphi^{(2q-k)},\varphi)_{L^2}\Big|\\
&\leq  \sum_{k=0}^{q}{q\choose k} {2p\choose
k}k!||x^{2p -k}\varphi||_2||\varphi^{(2q-k)}||_2\\
&\leq C  \sum_{k=0}^{q}{q\choose k} {2p \choose k}k!
(2/\theta)^{(q+p-k)}\frac{M^2_k}{M^2_k}M_{2q-k}M_{2p-k}\\
& \leq C \sum_{k=0}^{q}{q \choose k} {2p \choose k}
(2/\theta)^{(q+p-k)}\frac{k!}{M^2_k}M_{2q}M_{2p}\\
& \leq C H^{2(q+p)}M^2_{q}M^2_{p} \sum_{k=0}^{q}{q \choose k} {2p
\choose
k}(2/\theta)^{(q+p-k)}\\
&\leq C m^{q+p}M^2_{q}M^2_{p},
\end{split}
\end{equation}
\end{proof}

Let us now prove Theorem \ref{ljilja}.
\begin{proof}
(Theorem 2.2) (b) Let us suppose that
 for some
$\theta >0$
\begin{equation}
\sum_{n =0}^{\infty}|b_n|^2\exp[-2M(\theta\sqrt{2n+1})] < \infty.
\end{equation}
 We will prove the convergence of the series
 $\sum_{n=1}^{\infty}  \langle b_n h_n,\varphi \rangle $ for every $\varphi \in {\cal
 S}^{(M_p)}$. From that it follows that the mapping $f \mapsto
 \sum_{n=1}^{\infty}a_n b_n$ is an element of ${\cal S'}^{(M_p)}$.
 Using Schwartz inequality  we have
\begin{equation*}
\begin{split}
& \sum_{n=0}^{\infty}| \langle  b_nh_n,\varphi \rangle | =
\sum_{n=0}^{\infty}|b_n  \langle  h_n,\varphi \rangle |= \\
& =
\sum_{n=0}^{\infty}(|b_n|\exp[-M(\theta\sqrt{2n+1})]\cdot|  \langle  \varphi,h_n \rangle |\exp[M(\theta\sqrt{2n+1})])\\
& \leq \Big(\sum_{n=1}^{\infty}|b_n|^2\exp[-2M(\theta\sqrt{2n+1})]
\Big)^{1/2}\Big(\sum_{n=0}^{\infty}|  \langle \varphi,h_n \rangle
|^2\exp[2M(\theta\sqrt{2n+1})] \Big )^{1/2}.
\end{split}
\end{equation*}
The first  sum on the right hand side converges by supposition and
the convergence of the second sum follows from Theorem \ref{th 5}
since $\varphi \in{\cal S}^{(M_p)}$.

(a) Suppose that $f \in{{\cal S}^{(M_p)}}'$ and let $b_n =
  \langle  f,h_n \rangle $. Then for every  $\varphi = \sum a_nh_n \in{\cal
S}^{(M_p)}$, we have that the series $\sum a_nb_n$ converges.  We
will prove that for some $\theta > 0$
\begin{equation}
\sum_{n \in {\Bbb N}_0}|b_n|^2\exp[-2M(\theta\sqrt{2n+1})] <
\infty.
\end{equation}
Let us first prove that  the sequence
\begin{equation}\label{seq}
\{|b_n|\exp[-M(\theta\sqrt{2n+1})]\}_{n \in {\Bbb N}}
\end{equation}
is bounded for some $\theta$, say $\theta_0$. If it is not so,
then for every $q \in {\Bbb N}$ there exist $n_q \geq q$ such that
$$
|b_{n_q}|\exp[-M(q\sqrt{2n_q+1})] \geq 1.
$$
Let $\varphi = \sum a_mb_m$ where
$$
|a_m | =\exp[-M(q\sqrt{2n_q+1})]q^{-1}
$$
for $a_m= n_q$, and $a_m = 0$ for $m \neq n_q$, and $a_m$
satisfies that $a_mb_m = |a_mb_m|$. We will prove that $\varphi\in
{\cal S}^{(M_p)}$. Using \cite[Proposition 3.4, p.50]{K},
\begin{equation}\label{nejednakost}
M(k\rho) - M(\rho) \geq \frac{\log(\rho/A)\log k}{\log H}\geq
\frac{\log(\rho/A)\log k}{\log (1+H)}
\end{equation}
for $\rho=\mu \sqrt{2n_q+1}$, $k = q/\mu$, we have that
\begin{equation*}
\begin{split}
 \sum_{m=1}^{\infty}|a_m|^2& \exp[2M\mu\sqrt{2n+1}]\\
& = \sum_{q} q^{-2}\exp[-2M(q\sqrt{2n_q+1})+2M(\mu\sqrt{2n_q+1})]\\
& \leq \sum_{q}
q^{-2}\exp[-2\frac{\log(\frac{\mu}{A}\sqrt{2n_q+1})\log(\frac{q}{\mu})}{\log(1+H)}].
\end{split}
\end{equation*}
The above sum converges since for large enough $q$
$$
\exp[-2\frac{\log(\frac{\mu}{A}\sqrt{2n_q+1})\log(\frac{q}{\mu})}{\log(1+H)}]
$$
is bounded. From Theorem \ref{th 5} it follows that $\varphi \in
{\cal
 S}^{(M_p)}$, and therefore
 $$
 \sum_{n=0}^{\infty}a_nb_n =   \langle  f,\varphi \rangle  < \infty,
 $$
 which is in contradiction with
$$
 \sum_{n=0}^{\infty}a_nb_n =  \sum_{n=0}^{\infty}|a_nb_n| \geq
 \sum_{q=1}^{\infty}q^{-1} = \infty.
 $$
 Therefore, sequence (\ref{seq}) is bounded for some $\theta_0$.

 Using inequality (\ref{nejednakost}) for
 $\rho=\theta_0\sqrt{2n+1}$ and $k = \theta/\theta_0$, we have
 that
$$
\exp[-M(\theta\sqrt{2n+1})]\leq
\exp[-M(\theta_0\sqrt{2n+1})]\exp[-\frac{\log(\frac{\theta_0\sqrt{2n+1}}{A})\log(\frac{\theta}{\theta_0})}{\log(1+H}]
$$
which implies that the sequence
$$
b_n\exp[-M(\theta\sqrt{2n+1})]\rightarrow 0, \quad  n \rightarrow
\infty.
$$
for some $\theta > \theta_0$.

And finally, we will prove that the series (\ref{zaga}) converges
for some $\theta > \theta_0$. If we suppose that the sequence
diverges for every integer $\theta=q> \theta_0$, then there exists
an increasing sequence of integers $k_q$ suck that
$$
1  \leq   \sum_{n=k_{q-1}}^{k_q -1} |b_n| \exp[-M(q\sqrt{2n+1})]
<2 , \quad q = \theta_0+1, \theta_0+2,...
$$
 Put
$$
|a_n| = |b_n| \exp[-2M(q\sqrt{2n+1})]q^{-1}
$$
for $k_{q-1} \leq n,k_q$, $q  \rangle  \theta_0$. Then for every
fixed $\theta >0$ we have
\begin{equation*}
\begin{split}
&\sum_{n=k_{q-1}}^{k_q -1}|a_n|^2 \exp[2M(\theta\sqrt{2n+1})]=\\
& \sum_{n=k_{q-1}}^{k_q -1} |b_n|
\exp[-2M(q\sqrt{2n+1})]\exp[2M(\theta\sqrt{2m+1})-2M(q\sqrt{2n+1})].
\end{split}
\end{equation*}
Using once again the inequality (\ref{nejednakost})  we see that
$$
\exp[-2M(q\sqrt{2n+1})+2M(\theta\sqrt{2n+1})] \leq 1,
$$
and therefore, for every $\theta >0$,
\begin{equation*}
\begin{split}
\sum_{n=k_{q-1}}^{k_q -1}|a_n|^2 \exp[2M(\theta\sqrt{2n+1})]\leq
2q^{-2} .
\end{split}
\end{equation*}
It implies that the sum
$$
\sum_{n=0}^{\infty} |a_n|\exp[2M(\theta\sqrt{2n+1})]
$$
converges for every $\theta>0$. From theorem \ref{th 5}
 it follows that $\varphi=\sum a_n h_n$ belongs to ${\cal S}^{(M_p)}$
 which contradict the fact that the sum  $\sum a_nb_n$ diverges
 since
 $$
  \sum_{n=k_{q-1}}^{k_q-1} |a_nb_n| =
  \sum_{n=k_{q-1}}^{k_q-1}|b_n|^2\exp[-2M(q\sqrt{2n+1})]q^{-1}\geq q^{-1}.
   $$
\end{proof}

\bigskip
\address{
Department  of Mathematics and\\ Informatics\\
 University of Novi Sad \\
21000 Novi Sad\\
Serbia and Montenegro } {zlc@im.ns.ac.yu}
\address{
Department  of Mathematics and \\Informatics\\
 University of Novi Sad \\
21000 Novi Sad\\
Serbia and Montenegro } {dusanka@im.ns.ac.yu}
\end{document}